\documentclass{article}

\usepackage{arxiv}

\usepackage[utf8]{inputenc} 
\usepackage[T1]{fontenc}    
\usepackage{hyperref}       
\usepackage{url}            
\usepackage{booktabs}       
\usepackage{amsfonts}       
\usepackage{nicefrac}       
\usepackage{microtype}      

\usepackage{array}
\usepackage{hhline}
\usepackage{tabularx}
\usepackage{pdflscape}
\usepackage{geometry}

\usepackage{setspace}

\usepackage{enumitem}
\setlist{leftmargin=0.5in}
\usepackage{amsthm}
\usepackage{comment}
\usepackage{subcaption}

\usepackage{mathtools}
\DeclarePairedDelimiter{\ceil}{\lceil}{\rceil}

\usepackage{float}
\usepackage{multirow}
\usepackage{amsmath}
\usepackage{amsthm}
\usepackage{wasysym}
\usepackage{caption}
\usepackage{makecell}

\usepackage{enumitem}

\usepackage{algorithm}
\usepackage{algpseudocode}

\usepackage{framed} 
\usepackage{multicol} 

\title{Leveraging Owners' Flexibility in Smart Charge/Discharge Scheduling of Electric Vehicles to Support Renewable Energy Integration}

\author{
  Pouya ~Sharifi \\
  Department of Industrial \& Systems Engineering  \\
  Texas A\&M University\\
  College Station, TX \\
  \texttt{pouyasharifi@tamu.edu} \\
   \And
 Amarnath ~Banerjee \\
  Department of Industrial \& Systems Engineering\\
  Texas A\&M University\\
  College Station, TX\\
  \texttt{banerjee@tamu.edu} \\
    \AND
   Mohammad J.~Feizollahi \\
   Robinson College of Business \\
   Georgia State University, Atlanta, GA \\
   \texttt{mfeizollahi@gsu.edu} \\
}

\begin{document}
\maketitle

\begin{abstract}
High integration of intermittent renewable energy sources (RES), in particular wind power, has created complexities in power system operations. On the other hand, large fleets of Electric Vehicles (EVs) are expected to have great impact on electricity consumption, and uncoordinated charging process will add load uncertainty and further complicate the grid scheduling. In this paper, we propose a smart charging approach that uses the flexibility of EV owners to absorb the fluctuations in the output of RES in a vehicle-to-grid (V2G) setup. We propose an optimal scheduling algorithm for charge/discharge of aggregated EV fleets to maximize the integration of wind generation as well as minimize the charging cost for EV owners. Challenges for people participation in V2G, such as battery degradation and feeling insecure for unexpected events, are also addressed.
We first formulate a static model using mixed-integer quadratic programming (MIQP) with multi-objective optimization assuming that every parameter of the model is known a day ahead of scheduling. Subsequently, we formulate a dynamic (dis)charging schedule after EVs arrive into the system with updated information about EV availabilities, wind generation forecast, and energy price in real-time, using rolling-horizon optimization. Simulations using a group of 100 EVs in a micro-grid with wind as primary resource demonstrate significant increase in wind utilization and reduction in charging cost compared to uncontrolled charging scenario.
\end{abstract}


\section*{Nomenclature}
\noindent$T$ : Set of time/periods in scheduling with index $t$.\\
$\Delta t$ : Length of decision intervals.\\
$I$ : Set of all EVs with index $i$.\\
$I_{g2v}$ : Set of vehicles participating in G2V only.\\
$I_{v2g}$ : Set of vehicles participating in V2G.\\
$B$ : Set of vehicles that arrive with charge level below $SOC_{min}$.\\
$EV_{all}^t$ : Set of all vehicles plugged-in during time slot $[t,t+1]$.\\
$EV_{v2g}^t$ : Set of all vehicles in V2G plugged-in during time slot $[t,t+1]$.\\
$t_i^{arr}$ : Arrival time of EV $i$.\\
$t_i^{dep}$ : Departure time of EV $i$.\\
$T_i^{p}$ : Plug-in period of EV $i$.\\
$pr^t$ : Electricity price in time slot $[t,t+1]$ in $\cent$/kWh.\\
$D^t$ : Total charging demand in time slot $[t,t+1]$.\\
$W^t$ : Actual wind (renewable) production in time slot $[t,t+1]$ in kWh.\\
$W_{forecast}^t$ : Wind forecast for time slot $[t,t+1]$ in kWh.\\
$P_i^{c}$ : Maximum energy EV $i$ can take in $\Delta t$ period (kWh).\\
$P_i^{d}$ : Maximum energy EV $i$ can discharge in $\Delta t$ period (kWh).\\
$CP_i$ : Maximum charging power (kW) of the charger that EV $i$ is connected to.\\
$AR_i$ : Acceptance rate of EV $i$ in kW.\\
$\eta_i^{c}$ : Charging efficiency of EV $i$.\\
$\eta_i^{d}$ : Discharging efficiency of EV $i$.\\
$SOC_i^t$ : State of battery charge for EV $i$ at time $t$.\\
$SOC_{init,i}$ : Initial state of battery for EV $i$ in kWh.\\
$SOC_{cap,i}$ : Battery capacity of EV $i$ in kWh.\\
$SOC_{min,i}$ : Minimum level of battery charge for EV $i$ in kWh. \\
$SOC_{desired,i}$ : Desired level of battery charge for EV $i$ in kWh.\\
$T_{min,i}$ : Minimum number of periods to reach $SOC_{min,i}$. \\
$\Psi_i$ : The battery degradation cost for EV $i$ in $\cent$.\\
$\delta$ : Penalty for energy curtailment in $\cent$/kWh.\\
$\lambda$ : Owners' level of tolerance for battery degradation in [0, 1].\\
$P_G^{max}$ : Maximum transmission power between microgrid and the external grid in $\Delta t$ period.\\
$C_{bat, i}$ : Battery replacement cost in \$.\\
$X_{c,i}^t$ : Charging rate for EV $i$ in time slot $[t,t+1]$.\\
$X_{d,i}^t$ : Discharging rate for EV $i$ in time slot $[t,t+1]$.\\
$Y_{c,i}^t$ : Binary variable that takes a value of 0 if EV $i$ is charging in time slot $[t,t+1]$.\\
$Y_{d,i}^t$ : Binary variable that takes a value of 0 if EV $i$ is discharging in time slot $[t,t+1]$.\\
$Z_i$ : Binary auxiliary variable for EV $i$.\\
$\Omega^t$ : Wind curtailment in time slot $[t,t+1]$ in kWh.\\
$G^t$ : Energy supplied from the external grid in time slot $[t,t+1]$ in kWh.\\
$J$ : Set of planning times in dynamic modeling with index $j$.\\
$\Delta j$ : Planning intervals.\\
$A^j$ : Set of all vehicles arrived in time slot [$j-1$, $j$].\\
$E^j$ : Set of all vehicles to be planned at planning time $j$.\\
$E_{v2g}^j$ : Set of all vehicles in V2G mode to be planned at planning time $j$.\\
$E_{g2v}^j$ : Set of all vehicles in G2V mode to be planned at planning time $j$.\\
$\tau_{max}^j$ : The end of the planning (rolling) window when planned at time $j$.\\
$LC_i^{\{j\}}$ : The charging rate in the period previous to the planning time $j$.\\
$LD_i^{\{j\}}$ : The discharging rate in the period previous to the planning time $j$.\\
$E[ ]$: Expected value function.\\
$R_{v2g}$ : Ratio of vehicles participating in V2G.\\
$D_f^t$ : Charge demand for unknown future arrivals in time slot $[t,t+1]$.\\
$\hat{N_j^t}$ : Estimated number of vehicles arrived after time $j$ and plugged-in at time $t$.\\
$\hat{ER}$ : Estimated charge required during plug-in period.\\
$\Hat{PT}$ : Estimated plug-in period.\\
$P_{a,b}^k$ : The transition probability from wind state $a$ to $b$ in $k$ time steps.

\section{Introduction}\label{intro}

In recent years, more and more countries around the world are committing to reduce their carbon emissions and mitigate global warming. An integral part of cutting greenhouse emission is to reduce the burning of fossil fuels. To this end, RES are vital to replace or reduce the dependence on the existing gas and coal power plants. However, RES generation often depends on the weather, e.g. the sun for solar, wind for wind farms, rainfall for hydro, etc. While the grid can absorb a portion of RES generation, complete dependence on these intermittent sources requires massive grid-scale energy storage to cope with the intermittent nature of renewable generation. In addition to intermittent behavior of renewable electricity generators, their high initial cost is being a drawback. Thus, in many countries it has been actively promoted by the support of governments. Recent advances in technology are reducing their production cost and making them more competitive with conventional utilities.

As another solution to cut emission, countries are attempting to electrify the transportation sector as it accounts for large amount of Carbon emissions (20-25\% share) \cite{avci2014electric:2}. Electric vehicles look more promising than ever to replace the traditional internal combustion engine vehicles in the future since they could achieve zero emission if the electricity used is generated from renewable sources. However, the increasing adoption of EVs may cause a potential problem for the electric grid because of the unpredictable charging schedules of their owners. If unplanned, charging EVs will endanger the power system reliability and operations. The authors in \cite{yan2012impact:1'} examine the effects of EVs on power transformers in a distribution network. The study concludes that high penetration of EV will decrease the lifespan of the transformers significantly, as EV load increases the peak load and that can exceed the capacity of transformers. It suggests that charging EVs during nighttime can decrease the transformers loss-of-life factor. The authors in \cite{verzijlbergh2011impact:10'} compare the uncontrolled and controlled charging scenarios and conclude that the proposed controlled scenario can reduce transformers' overload by 25\% with a large penetration (max of 75\% in 2040) of EVs. The increase in EV penetration leads to not only an increase in electricity demand but also the shape of the demand will change significantly resulting in higher demand variability and impacting electricity infrastructure and making it difficult to accurately predict the load. Peak demand determines the system capacity requirements; thus, increasing the peak demand will affect the electricity infrastructure of power system \cite{muratori2018impact:7ev}. Most importantly, uncoordinated charging of EVs adds a stochastic element to the power system, which complicates the planning and operation of power systems (e.g., unit commitment and economic dispatch).

In anticipating this problem, many have proposed the concept of vehicle-to-grid (V2G) to integrate these future EVs into the grid successfully. V2G sees the EVs as a distributed generation/storage system as well as dynamic flexible load which could be utilized to balance the supply and demand of electricity \cite{mwasilu2014electric:15ev}. V2G technology can improve the power system operations along several services, such as frequency/voltage regulation, spinning reserve, and peak shaving \cite{tan2016integration:12}. The work in \cite{sheikhi2013strategic:33} presents a game-theoretic approach in load management strategy for EV charging to reduce peak load considering dynamic behavior of EV drivers as well as electricity price. However, intermittent generation of RES is not considered in this work.

With the focus on maximizing RES integration, some researchers have attempted to find the optimal scheduling for EV charging. In \cite{borba2012plug:14ev}, the authors formulate the optimal charging schedule in charge-only mode and find that charging EVs overnight can absorb the excess power generated by wind power; thus, increasing the RES utilization. In \cite{jin2013energy:14'}, the authors present an optimization algorithm for charging schedule problem to minimize energy cost from the grid based on a queuing model that does not require any information or prediction about the wind production, EV charging request, and electricity price. The proposed approach leads to 78 \% cost reduction compared to the "charge-upon-arrival" case. The study in \cite{schuller2015quantifying:34ev} formulates the EV charging schedule as a mixed integer programming (MIP) problem maximizing the RES integration for the day-ahead problem assuming perfect generation forecast for the planning horizon. Most articles in this domain have considered the bidirectional V2G in which the EVs are also capable of providing discharge. 
In \cite{gottwalt2013assessing:17'}, the authors developed a MIP problem to maximize the utilization of renewable energy while satisfying EV and household demand assuming known future trips as well as EV and household loads in day-ahead scheduling framework. Likewise, the authors in \cite{ghofrani2012electric:28ev} propose an energy management system to maximize the utilization of wind energy assuming known wind and load. In this work, the storage capability of EVs in the context of a distributed feeder with primary wind resource is examined. It is shown that wind utilization can reach 89\% with a coordinating charging and discharging strategy under the assumption of an EV for each household. 
In \cite{honarmand2014optimal:13ev}, the authors propose an intelligent scheduling algorithm with the objective of maximizing owners' profit under the assumption of known day-ahead arrival and departure times. The work in \cite{lopes2009using:11ev} proposes to utilize the storage capability of a large fleet of EVs that can be used as distributed storage units to help keep the grid frequency within a certain limit in the presence of intermittent RES generation in a distributed grid framework. 

The EV arrival and departure time are dynamic and time-varying since the driving behavior of EV owners is uncertain as it depends on several factors, such as traffic conditions, and randomness in commuting behavior of drivers \cite{wan2018model:31'}. In this work, we propose a real-time multi-objective optimization framework in which the EV characteristics and behavior are only known upon vehicle arrivals, not a day ahead of scheduling. We developed a dynamic planning horizon algorithm for the aggregator for EV charging and discharging. Our work is similar to \cite{he2012optimal:10ev}, where the authors have considered the dynamic behavior of EV owners in an optimal charge/discharge scheduling of EVs using rolling horizon optimization method; however, the intermittent generation of RES is not considered in their work. To evaluate our dynamic scheduling, we first model the optimal schedule for a day-ahead static scenario, where commuting schedule of all EVs are assumed to be known day ahead of time. The results of the dynamic model is then compared to the static model. Different EV characteristics and flexibility of owners for different wind and electricity price profiles are studied to evaluate the performance of our model.

Moreover, battery degradation has not been taken into consideration in any of the mentioned articles (except in \cite{he2012optimal:10ev}). There are three main challenges in V2G that need to be addressed before applying it to the real-world:
\begin{enumerate}[topsep=0pt,itemsep=-1ex,partopsep=0ex,parsep=1.5ex]
\item \textit{Battery degradation}:
One major impediment to V2G is its toll on the battery life \cite{peterson2010economics:23'}. An average lithium ion battery utilized by major automobile companies lasts between 2000-3000 charging cycles \cite{castillo2014grid:86}, cycling them daily will significantly reduce their lifespan, hence the economic benefits might be completely countered by the battery wear.
\item \textit{Efficiency}:
V2G might not be the most efficient method to store electricity as well. Round trip power losses of all EV components are reported to be between (17-25\%), which means V2G as an energy storage system can only return 75-83\% of its energy input \cite{apostolaki2017measurement:29'}. Comparing to other ways of energy storage (e.g., pumped hydro storage 85-95\%), V2G as an energy storage system might lose more in transaction.
\item \textit{Feeling insecure for urgent needs}:
The owner of an EV may not be willing for a third party to directly control their charging process because in case of emergency or unexpected travel need the battery might have insufficient range to carry the vehicle to its destination \cite{tan2016integration:12}. 
\end{enumerate}

In this study, we address all these challenges in V2G and people participation that discourage people to allow a third party to control their charging process. Unlike some of the previous research articles, the proposed algorithm does not attempt to start or stop charging EVs frequently, which leads to battery degradation. The approach here is based on optimal scheduling of a large number of EVs depending on the availability of EVs and preferences, needs and flexibility of their owners.

We study the EV scheduling for both static and dynamic models in two cases: 1) "Business-as-usual" (BAU) case: immediate charge upon arrival at full speed until reaching full charge; 2) smart schedule that an aggregator can modify the for charge/discharge rate in discrete time intervals. We also consider different ratio of EV owners that participate in V2G to evaluate the benefits of V2G compared to G2V.

The remainder of this paper is structured as follows. In Section \ref{SystemModel}, we provide an overview of our system model and the proposed approach. In Section \ref{sec:Static}, we formulate the optimal scheduling problem of static model, while the dynamic modeling is provided in Section \ref{Dynamic}. Section \ref{sec:simulation} discusses the simulation results of our study. Finally, we conclude the paper in Section \ref{conclusion}.

\section{Approach and System Model Overview} \label{SystemModel}

\begin{figure}[t!]
    \centering
    \caption{Smart grid communciations}
    \includegraphics[scale=0.75]{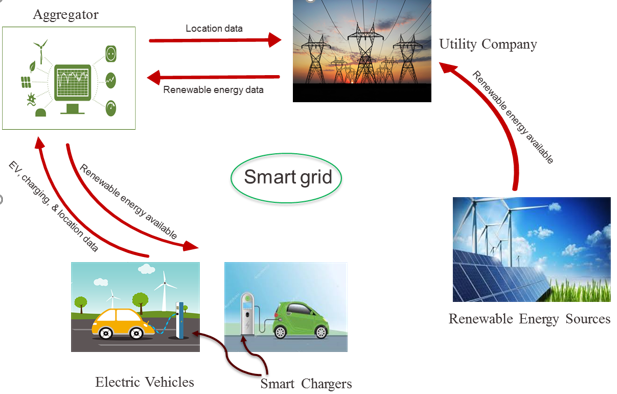}
    \label{fig:smartgrid}
\end{figure}

Smart energy metering and advanced controls have enabled real-time communication in smart grid \cite{mwasilu2014electric:15ev}. We consider a smart grid with real-time communication between its participants, and an aggregator who is responsible to offset the fluctuation of RES (e.g., wind energy) with the optimal charging/discharging schedule of EVs. A schematic of the relationship is shown in Figure \ref{fig:smartgrid}. The aggregator, which can be a utility company or a third party, gets information about available renewable energy, and market electricity price. EVs are connected to the smart chargers, and the aggregator receives data related to the status of EVs, the charge requirements, and their owners' preferences. The aggregator is responsible for selling power output of the renewable sources as much as possible with a price less than what conventional utilities offer. On the other hand, suppose that there is a set of EV owners who want to charge their vehicles with renewable energy as much as possible for sustainability concerns. They may also have an economic motivation to pay less for charging their vehicles. Since EVs are idle 90\% of the time but they require a few hours to recharge, the goal is to use the flexibility of EVs and their drivers to absorb the fluctuations in the output of RES, while satisfying their requirements and concerns. The hypothesis is that the EVs charge their batteries with energy from RES and discharge energy in low wind periods. Also, in case where there is insufficient wind to meet the demand during the whole plug-in period, the EV charging process should be shifted to low electricity price periods whenever possible.

\subsection{Electricity Generation and Consumption Model}
We study a MicroGrid (MG) with wind power as primary resource and local parking lots as consumers. EV charging in parking lots are the only demand for wind generation and there is no storage unit in the system, meaning that if wind production is more than the total charging demand, excessive wind energy has to be curtailed. Since wind energy has near zero marginal generation cost and because of the support policy from the government, we assume that wind energy has no cost. The microgrid is connected to the external grid via transformer for back-up power, so that if charging demand exceeds the available wind energy, the remaining energy is purchased from the grid with the real-time electricity price. EV owners can participate in either unidirectional or bidirectional V2G. The former, also referred to as grid-to-vehicle (G2V), is when the EV can only charge energy from the electrical grid. However, in bidirectional V2G, the EVs can inject energy to the gird by providing discharge. In the context of this problem, V2G is referred to bidirectional flow of energy between energy source and vehicle. In case of low wind generation, EVs in V2G mode (set $I_{v2g}$) can provide energy by discharging their batteries to charge other EVs. We assume that the discharge energy is only used for charging other EVs, it is not sold to the external grid. The EV discharge energy is sold with a price slightly less than that of real-time wholesale market. By that, we make sure that EVs are making revenue from selling their stored energy and also other EVs are paying less compared to the market price. The maximum charging power for an EV in one hour interval is calculated as $ \min (AR_i,\;\; {CP}_i)$, where $AR_i$ is the maximum power the EV can take, and $CP_i$ is the maximum charging power of the outlet that EV $i$ is connected to. Thus, the maximum energy EV $i$ can take in each $\Delta t$ decision period is defined by $P_i^{c}=\min (AR_i,\;\; CP_i)\times (\Delta t/1 hour)$. We assume that charging and discharging powers are the same ($P_i^{c}=P_i^{d}$).

Local consumers in the MG are the EV owners parking their vehicles in residential places or in workplace parking lots. For a total of N vehicles, we assume half of them parked at workplace and the other half are the vehicles charging their EVs at home, and vehicles are plugged-in as soon as they arrive. We have excluded the consumption for households, since in case of known or close to known household demand, that would not add value to our optimization problem.

\subsection{Battery Degradation Model} \label{battDeg}
Researchers in \cite{peterson2010economics:23'} have found that the maximum annual profit for an EV is very limited if considering the battery degradation cost, and without considering this cost, the profit is exaggerated. Moreover, battery degradation contributes as a great challenge in people participation, so we study the battery degradation cost in the charge/discharge process to reduce the total cost for the owners. Finding exact degradation cost for any battery is out of scope for this paper. Here, we consider two models for battery degradation used in literature. The first model considers a quadratic function, which consists of two terms, one for charge/discharge rate and the other term captures the cost of degradation for fluctuation in energy rate \cite{he2012optimal:10ev}. The second model is adopted from \cite{peterson2010economics:23',sortomme2011optimal:11'}, which modeled the battery wear cost as a linear function of battery replacement cost and percent of battery used. The laboratory measurements in \cite{peterson2010economics:23'} predicted a cost of 4.2 $\cent/kWh$ for a battery pack with \$5,000 replacement cost. Thus, in this model, we consider the same value for degradation cost coefficient. Equations (\ref{eq:1000})-(\ref{eq:1001}) define the degradation cost of EV $i$ in the second model.
\begin{eqnarray}
 \Psi_i = \sum_{t \in T_i^p} 4.2\times \frac{ C_{bat, i}}{5,000} \times (percent \; of \; battery \; used)^t \label{eq:1000}\\
 (percent \; of \; battery \; used)^t = \frac{P_{c,i}X_{c,i}^t-P_{d,i}X_{d,i}^t}{SOC{cap,i}} \label{eq:1001}
\end{eqnarray}

Our emphasis is on minimizing the cost of battery degradation caused by frequent recharge and discharges as well as the energy processed by the EV, so we mainly use the first model (presented in equations (\ref{eq:5})-(\ref{eq:6})) for our reference, but the second model is also used for validation of results. We penalize the battery degradation cost in the objective function so that the optimal solution reduces the frequency of stop and start of the charge process.

\section{Static Scheduling Optimization} \label{sec:Static}
In this section, a scheduling algorithm is proposed for the aggregator that determines the day ahead EV (dis)charge schedule. In the static model, the following assumptions are made:
\begin{itemize}
  \setlength\itemsep{0.0em}
  \item The EV owners are obligated to provide their arrival time, departure time, the desired level of charge, and minimum level of charge, the night before of the scheduling day.
  \item The initial state of charge and the EV characteristics, such as battery capacity, acceptance rate, and battery replacement cost are known.
  \item The wind production and electrical grid energy price for the next day is forecasted with perfect accuracy.
\end{itemize}
Given the above information, the aggregator determines the charging schedule for the next planning horizon (e.g., next day). We consider the scheduling  for discrete time intervals of equal size $\Delta t$. This smart scheduling is determined to maximize wind utilization, minimize the grid supply, and minimize the charging cost for EV owners. Comparison between the smart charging case versus BAU case, where the vehicles are charged immediately upon their arrivals until full charge, is provided.
Although the assumptions made in this deterministic case are unrealistic, this model provides the global optimal solution for the case where all model parameters are known. The results of the static model will be used later in this paper to evaluate the performance of the dynamic model.

\subsection{Modeling \& Mathematical Formulation}

\begin{figure}[h!]
    \centering
    \caption{Illustration of difference between time and period}
    \includegraphics[scale=0.75]{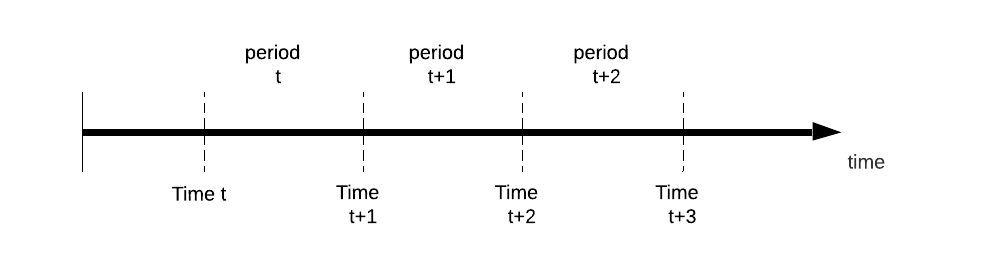}
    \label{fig:time&period}
\end{figure}

We study the EV charging and discharging during the planning horizon that is evenly divided into time intervals of $\Delta t$ minutes. The aggregator finds the optimal charge/discharge rate for all EVs in each period $t$. In this paper, the time slot [$t,t+1$] is referred to as period $t$ (see Figure \ref{fig:time&period}). It is assumed that grid electricity price, and charge/discharge rate remain constant for the entire interval of $\Delta t$ minutes. Wind power availability in $\Delta t$ time interval is known and is utilized to charge the EV fleet during that interval. Here, we use 15-minute interval, which divides a single day into 96 equal intervals.

The decision variables used in the model are described as follows. $X_{c,i}^t$ and $X_{d,i}^t$ are continuous variables between 0 and 1 that determine the charge and discharge rate of EV $i$ in period $t$ with a value of 1 meaning full-speed charging/discharging, and 0 meaning remaining untouched. The state of charge for EV $i$ at time $t$ is captured by $SOC_i^t$. $G^t$ and $\Omega^t$ are grid supply and wind curtailment in period $t$, respectively. Finally, $Z_i$, $Y_{c,i}^t$, and $Y_{d,i}^t$ are binary auxiliary variables.

The objective function is to minimize a linear combination of the charging cost of energy purchased from the electric grid, the battery degradation cost, and the wind curtailment penalty cost.
\begin{equation}
 \displaystyle{\min_{X_c, X_d} \quad \sum_{t \in T} pr^t G^t + \sum_{i \in I} \lambda \Psi_i + \sum_{t \in T} \delta \Omega^t} \label{obj:1}
\end{equation}

\noindent The variable $G^t$ is the energy supplied from the external grid (conventional generators). It can be defined by inequality (\ref{eq:2}) and non-negativity constraint (\ref{eq:3.24}).  
\begin{eqnarray}
 G^t \geq \displaystyle \sum_{i \in EV_{all}^t} (X_{c,i}^t P_i^{c}) - \displaystyle \sum_{i \in EV_{v2g}^t} (X_{d,i}^t P_i^{d}) - W^t, \quad \forall t \in T \label{eq:2}
\end{eqnarray}
It can be easily noted that $G^t$ is the maximum of zero and the total net charging demand minus the wind energy in any period $t$. Inherited in the first term of the objective function, the optimization problem attempts to minimize the charging cost for EV owners that is purchased from the grid. Also, imposing charging cost for the grid energy causes the demand to seek cheap wind energy and increases wind utilization.

The next term is the total battery degradation cost multiplied by $\lambda$, which is a hyper-parameter that indicates the level of tolerance of the owners for their battery degradation. Value of $\lambda=1$ means no tolerance for the owner. Battery degradation is calculated using a quadratic function as defined in (\ref{eq:5})-(\ref{eq:6}) \cite{he2012optimal:10ev}.
\begin{eqnarray}
 \Psi_i = \sum_{t \in T_i^p} \alpha(\eta_i^{c} P_i^c (X_{c,i}^t  -X_{c,i}^{t-1}))^2 + \beta(\eta_i^{c} P_i^c X_i^t)^2, \quad \forall i \in I_{g2v}  \label{eq:5}\\
 \begin{split}
 \Psi_i {}={} \sum_{t \in T_i^p} \alpha[\eta_i^{c} P_i^{c}  (X_{c,i}^t -X_{c,i}^{t-1})]^2 + \beta[\eta_i^{c} P_i^{c}  X_{c,i}^t]^2 \\+
 \alpha[P_i^{d}/\eta_i^{d} (X_{d,i}^t -X_{d,i}^{t-1})]^2 + \beta[(P_i^{d} / \eta_i^{d}) X_{d,i}^t]^2, \quad \forall i \in I_{v2g} \label{eq:6}\\
\end{split}
\end{eqnarray}
 

In the third term, a small penalty $\delta$ is considered for wind curtailment. The wind curtailment denoted by $\Omega^t$ is the maximum of zero and wind production minus the charging demand at each period. $\Omega^t= \max\{0, W^t - D^t\}$ (see constraint (\ref{eq:4})).
\begin{eqnarray}
 & \Omega^t \geq W^t + \displaystyle \sum_{i \in EV_{v2g}^t}(X_{d,i}^t P_i^{d}) -\displaystyle \sum_{i \in EV_{all}^t}(X_{c,i}^t P_i^{c}), \quad \forall t \in T \label{eq:4}
\end{eqnarray}

Note that the demand for wind energy in period $t$ can be easily calculated by $W^t-\Omega^t$. It is good to mention that in each time interval $t$, at most one of $\Omega^t$ and $G^t$ can be positive. Since the EV charge load is the only demand for wind energy, minimizing wind curtailment implies maximizing wind utilization. Similar to $\lambda$, $\delta$ is a hyper-parameter that determines the weight for wind curtailment penalty. Considering a small value for $\delta$, this term comes into play only when there is enough wind production, and it ensures to reduce the wind curtailment by charging the vehicle to their full battery capacity instead of the user specified desired level.

We assume a linear charging behavior for the batteries. Thus, the state of charge is initialized and updated in a set of constraints (\ref{eq:3.1})-(\ref{eq:3.3}).
\begin{eqnarray}
 & SOC_i^{t_i^{arr}} = SOC_{init,i} & \forall i \in I \label{eq:3.1}\\
 & SOC_i^t = SOC_i^{t-1} + \eta_i^{c} P_i^{c} X_{c,i}^{t-1} &\forall i \in I_{g2v}, \;\; \forall t \in T_i^p \label{eq:3.2}\\
 & SOC_i^t = SOC_i^{t-1} + \eta_i^{c} P_i^{c} X_{c,i}^{t-1} - P_i^{d} X_{d,i}^{t-1} /\eta_i^{d}   &\forall i \in I_{v2g}, \;\; \forall t \in T_i^p \label{eq:3.3}
\end{eqnarray}

Constraint (\ref{eq:3.1}) sets the initial state of battery charge upon arrival to the state of charge at arrival time ($SOC_i^{t_i^{arr}}$). $SOC$ at each time is updated in constraints (\ref{eq:3.2})-(\ref{eq:3.3}) by adding the charging energy for the vehicle in the current period to the charge level of the previous time. Constraint (\ref{eq:3.7}) limits the total grid supply by the transformer's capacity.
\begin{eqnarray}
 & G^t \leq P_{max}^{G} & \forall t \in T \label{eq:3.7}
\end{eqnarray}

The next two constraints (\ref{eq:3.9})-(\ref{eq:3.10}) guarantee that if the vehicle cannot reach the desired level in its designated charging period, the vehicle is charged with full speed the entire plug-in period. However, if the vehicle can reach the desired level, constraint (\ref{eq:3.11}) makes sure that the final state of charge (SOC at departure time) is at least as much as the desired level of battery charge requested by the user.
\begin{eqnarray}
 & SOC_{init,i} + P_i^{c}\eta_i^{c} (t_i^{dep}-t_i^{arr}) \geq SOC_{desired,i} - MZ_i & \forall i \in I \label{eq:3.9}\\
 & X_{c,i}^t \geq Z_i & \forall i \in I, \;\; \forall t \in T_i^p \label{eq:3.10}\\
 & SOC_i^{t_i^{dep}} \geq SOC_{desired,i} - M Z_i & \forall i \in I \label{eq:3.11}
\end{eqnarray}

Constraint (\ref{eq:3.12}) limits the state of charge at any time to the battery capacity of an EV.
\begin{eqnarray}
 & SOC_i^t \leq SOC_{cap,i} & \forall i \in I, \;\; \forall t \in \{t_i^{arr},..., t_i^{dep}\} \label{eq:3.12}
\end{eqnarray}

The next group of constraints sets a minimum level of charge for all vehicles. If the vehicles arrive with a charge level less than the minimum level (set $B$), then it has to charge with full speed to get to the minimum level as soon as possible (Eq. (\ref{eq:3.14})). After reaching the minimum state of charge, the SOC should never drop below $SOC_{min}$ (constraint (\ref{eq:3.15})).
\begin{eqnarray}
 & SOC_i^t \geq SOC_{min,i} & \forall i \in I \setminus B, \;\; \forall t \in \{t_i^{arr},..., t_i^{dep}\} \label{eq:3.13}\\
 & X_{c,i}^t = 1 & \forall i \in B, \;\; \forall t \in \{t_i^{arr},..., t_i^{arr} + T_{min,i}\} \label{eq:3.14}\\
 & SOC_i^t \geq SOC_{min,i} & \forall i \in B, \;\; \forall t \in \{t_i^{arr} + T_{min,i},..., t_i^{dep}\} \label{eq:3.15}
\end{eqnarray}

$T_{min,i}$ denotes the minimum number of periods that the EV $i$ has to charge with full speed to reach $SOC_{min,i}$, and is calculated by equation (\ref{eq:Tmin}).
\begin{equation}
 T_{min, i}  =  \ceil[\bigg] {\frac{SOC_{min,i} - SOC_{init,i}}{\eta_i^c P_i^c}}
 \label{eq:Tmin}
\end{equation}
Constraints (\ref{eq:3.16})-(\ref{eq:3.16'}) set the charge/discharge rates to zero for the period prior to arrival.
\begin{eqnarray}
 & X_{c,i}^{t_i^{arr}-1} = 0 & \forall i \in I_{g2v} \label{eq:3.16}\\
 & X_{c,i}^{t_i^{arr}-1} + X_{d,i}^{t_i^{arr}-1}  = 0 & \forall i \in I_{v2g} \label{eq:3.16'}
\end{eqnarray}
Constraints (\ref{eq:3.17})-(\ref{eq:3.19}) ensure that in any period during plug-in time, the vehicles in V2G mode can either charge, discharge, or do nothing.
\begin{eqnarray}
 & X_{c,i}^{t} \leq 1-Y_{c,i}^t & \forall i \in I_{v2g}, \;\; \forall t \in T_i^p   \label{eq:3.17}\\
 & X_{d,i}^{t} \leq 1-Y_{d,i}^t & \forall i \in I_{v2g}, \;\; \forall t \in T_i^p   \label{eq:3.18}\\
 & Y_{c,i}^t +Y_{d,i}^t = 1 & \forall i \in I_{v2g}, \;\; \forall t \in T_i^p   \label{eq:3.19}
\end{eqnarray}
Finally, constraints (\ref{eq:3.20})-(\ref{eq:3.24}) specify the binary and non-negativity constraints for the decision variables. 
\begin{eqnarray}
 & 0 \leq X_{c,i}^t \leq 1  & \forall i \in I, \;\; \forall t \in T_i^p \label{eq:3.20}\\
 & 0 \leq X_{d,i}^t \leq 1 & \forall i \in I_{v2g}, \;\; \forall t \in T_i^p \label{eq:3.21}\\
 & Y_{c,i}^1,\; Z_i \;\; Binary & \forall i \in I, \;\; \forall t \in T_i^p\label{eq:3.22}\\ 
 & Y_{d,i}^t \;\; Binary & \forall i \in I_{v2g}, \;\; \forall t \in T_i^p \label{eq:3.23}\\ 
 & \Omega^t,\; G^t \geq0 & \forall t \in T \label{eq:3.24}
\end{eqnarray}

Although, the assumptions made in this deterministic case were somewhat impractical, this model provides the global optimal solution for the case where all model parameters are known. It also provides a baseline for the dynamic model which we will describe in the next section. In the dynamic model case, we address some of the impracticalities of the static model.

\section{Dynamic Scheduling Optimization} \label{Dynamic}
In a more realistic scenario, the assumption for the obligation of providing perfect information a day ahead by the EV owners is relaxed. In this dynamic model, EV owners input their needs and preferences (departure time, desired level of battery, minimum required level of charge) upon arrival. The smart chargers, on the other hand, automatically detect the necessary EV characteristics, such as battery capacity, acceptance rate, and state of charge. At any planning time $j$, the aggregator batches all the vehicles that have arrived during the $[j-1, j]$ time slot. The aggregator also considers those vehicles that have not departed from the previous periods and are still in charging. With updated information regarding the number of EVs and their requirements, the renewable energy generated, and the price of electricity, the scheduling algorithm (which is discussed in \ref{algorithm}) optimizes the charging schedule for the current planning window (rolling window). The planning window is defined by the period from time $j$ till the time that all vehicles in set $E^j$ departs. $E^j$ is the set of vehicles considered in planning at time $j$. This schedule is called dynamic because it can be updated as wind production forecast, the electric price forecast, and EV availabilities are updated. Also, note that the charging schedule of all vehicles gets updated frequently at each planning time until they depart, thus minimizing the effect of uncertainty in wind generation and electricity price. We will consider the uncertainty of wind generation forecast in \ref{sec:intermittentWind}. 

\subsection{Algorithm} \label{algorithm}

Considering 1-hour planning intervals $\Delta j$, one day is divided equally into 24 planning times so that planning times are at each exact hour of the day. For ease of notation, $j+1$ denotes the next planning time, which is one hour after current planning time $j$. However, to be consistent with the notation of our decision periods $t$, which increment by one every $\Delta t= 15$ minutes, another variable is defined by $\phi_j= 4j$. For instance, $j=0 \implies \phi_j= 0$ denotes planning at time 12:00 am, and $j=1 \implies \phi_j=4$ denotes the next planning time at 1:00 am.

\begin{figure}[t]
    \centering
    \caption{Illustration of planning (rolling) window in system with three arrivals}
    \centerline{\includegraphics[scale=0.65]{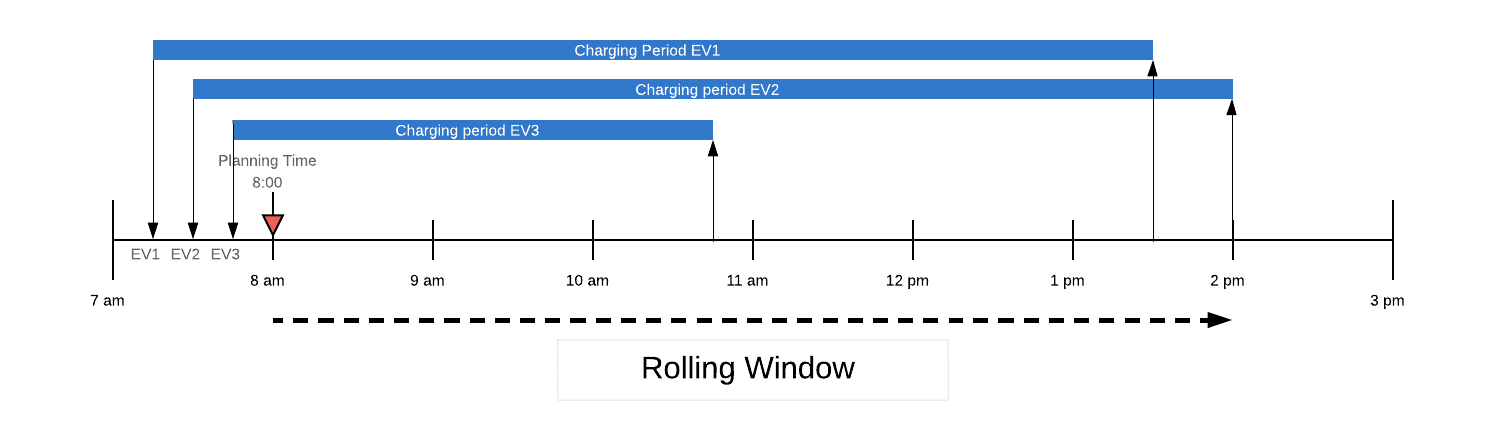}}
    \label{fig:threeArrivals}
\end{figure}

At any planning time $j$, we have a set of arrivals during time slot $[j-1, j]$, called $A^j$. For instance, $A_1$ consists of all the vehicles arriving between 12:00 am and 1:00 am. Vehicles in $A^j$ are added to the set $E^j$, which accounts for the vehicles that need to be planned at time $j$. The aggregator also batches those vehicles already arrived and planned in the previous period $j-1$ (vehicles in set $E_{j-1}$) if the vehicle remains plugged-in after time $j$. The vehicle is also added to the set $E^j$. Based on owners' input, the aggregator determines the planning window, which is from time $j$ to the time that last vehicle departs, denoted by $\tau_{max}^j$. It is defined by $\tau_{max}^j = \max\{t_i^{dep}|\;\; i\in E^j\}$. As an example, consider planning for EVs at $j = 8$ (8 am). In Figure \ref{fig:threeArrivals}, arrival and departure times for vehicles \{1,2,3\} are shown by down arrows and up arrows, respectively. Set $A_8 = \{1,2,3\}$ is the set of EVs arriving in period $j= $ [7 am, 8 am]. Assuming that there are no arrivals before 7 am, the set $E_8$ is equal to $A_8$. The charging (plug-in) period for EVs are depicted as blue lines. The end of planning window, $\tau_{max}^8$, is defined by $\tau_{max}^8 = \max\{t_i^{dep}|\;\; i\in E_8\}$. The last departure is when vehicle 2 departs at 2 pm ($t_2^{dep}=56$). Thus, the planning window is from $\phi_8 = 32$ to $\tau_{max}^8 = 56$.

If the arrivals after $j$ were known to us, the problem and optimal solution would be more accurate resulting in a global optimal solution. Since they are unknown, we need to estimate the average charge demand of future arrivals for the current planning window. To estimate the amount of charge, we need to calculate the average number of EVs in each $\Delta t$ period, the average charge required, and the average plug-in period by analyzing the historical data.
\begin{equation}
  E[D_f^t] = \frac{\hat{ER}\cdot \hat{N_j^t}}{\hat{PT}}
  \label{FutureDemand}
\end{equation}
\noindent $\hat{ER}$ is the estimated charge required by an EV, and it can be calculated by $\hat{ER} = \Hat{SOC}_{des}- \hat{SOC}_{init}$. After finding the set of vehicles to plan at time $j$, their state of charge, the charge/discharge rate for the period prior to $j$, and the estimated future charge demand, the aggregator runs the optimization algorithm that solves the MIQP problem (which is described in \ref{sub:modelUniDyn}) to find the optimal charging procedure for each EV. The aggregator can repeat this process for the next planning time $j+1$. This repetitive algorithm will run for the planning horizon $J$, which can be from a few hours to a couple of years. The steps of the rolling horizon approach are provided in Algorithm \ref{DynamicAlgorithm}.

\begin{algorithm}[t!]
\caption{Rolling Horizon Algorithm}
\label{DynamicAlgorithm} \begin{algorithmic} 
\State Initialize $E^j = \emptyset, \;\; \forall j \in J+\{0\}$, and set $j = 1$.
\For{$j \in J$}
\State $E^j \gets A^j$
\For{$i \in E_{j-1}$}
    \If{$t_i^{dep} > \phi_j$}
        \State Add $i$ to the $E^j$. $E^j \gets E^j +\{i\}$
    \EndIf
\EndFor
\For{$i \in E^j$}
        \State Set $LC_i^{\{j\}}, LD_i^{\{j\}}$ as the charge/discharge rate in the previous period.
        \State Update the $SOC_{init,i}$ as the $SOC$ of EV $i$ at time $j$.
\EndFor
\State Update the uncertain parameters ($W^t$ and $pr^t$) with their recent forecast.
\State Set $B = \{i|\;\; i \in E^j,\;\; SOC_{min,i} > SOC_{init,i}\}$
\State Find the end of current planning window by $\tau_{max}^j = \max\{t_i^{dep}|\;\; i\in E^j\}$
\State Calculate the expected charge demand of future arrivals from equation (\ref{FutureDemand}).
\State Run the optimization algorithm described in \ref{sub:modelUniDyn}.
\State Charge and discharge the EVs according to optimal $X_{c}, X_d$ values.
\EndFor
\end{algorithmic} 
\end{algorithm}

\subsection{Modeling \& Mathematical Formulation} \label{sub:modelUniDyn}
At each planning time $j$, the algorithm needs to solve a MIP or MIQP problem (depending on the battery degradation model) similar to the static case. The optimization problem at planning time $j$ considers the following objective function:
\begin{equation}
 \displaystyle{\min \quad \sum_{t=\phi_j}^{\tau_{max}^j-1} pr^t G^t + \sum_{i \in E^j} \lambda \Psi_i + \sum_{t=\phi_j}^{\tau_{max}^j-1} \delta \Omega^t} \label{obj:2}
\end{equation}

The objective function (\ref{obj:2}) is similar to the static case  (\ref{obj:1}), except for the fact that it minimizes the cost for the planning window, which is $\{\phi_j, ..., \tau_{max}^j-1\}$ instead of the entire planning horizon $T$. Here, $\tau_{max}^j$ is the end of planning window and is defined by $\tau_{max}^j = \max\{t_i^{dep}|\;\; i\in E^j\}$. Also, note that the total plug-in period ($T_i^p$) is no longer from arrival time to departure time, it is from planning time $j$ till departure. Thus, $T_i^p$ is updated by \{$\phi_j , ..., t_i^{dep}-1$\}.
\begin{eqnarray}
 \Psi_i = \sum_{t = \phi_j}^{t_i^{dep}-1} \alpha(\eta_i^{c} P_i^c  (X_{c,i}^t  -X_{c,i}^{t-1}))^2 + \beta(\eta_i^{c} P_i^c X_i^t)^2, \quad \forall i \in E^j_{g2v}  \label{eq:100}\\
\begin{split}
 \Psi_i = \sum_{t = \phi_j}^{t_i^{dep}-1} \alpha[\eta_i^{c} P_i^{c}  (X_{c,i}^t -X_{c,i}^{t-1})]^2 + \beta[\eta_i^{c} P_i^{c}  X_{c,i}^t]^2 +
 \alpha[P_i^{d}/\eta_i^{d} (X_{d,i}^t -X_{d,i}^{t-1})]^2 \\
  +\beta[(P_i^{d}/\eta_i^{d}) X_{d,i}^t]^2, \quad \forall i \in E^j_{v2g} \label{eq:101}
\end{split}\\
 G^t \geq \displaystyle \sum_{i \in EV_{all}^t} (P_i^{c} X_{c,i}^t) + E[D_f^t] - \displaystyle \sum_{i \in EV_{v2g}^t} (P_i^{d} X_{d,i}^t) - W_{forecast}^t, \quad \forall t \in \{\phi_j, ..., \tau_{max}^j-1\} \label{eq:102}\\
 \Omega^t \geq W_{forecast}^t + \displaystyle \sum_{i \in EV_{v2g}^t}(P_i^{d} X_{d,i}^t ) -\displaystyle \sum_{i \in EV_{all}^t}(P_i^{c} X_{c,i}^t) - E[D_f^t], \quad \forall t \in \{\phi_j, ..., \tau_{max}^j-1\} \label{eq:103}
\end{eqnarray}

The problem is subject to a set of constraints, most of which are similar to the static case with a few modifications. Since planning occurs at time $j$, only vehicles in the set $E^j$ from time $\phi_j$ to $\tau_{max}^j$ are included in the problem. Constraints (\ref{eq:2.1})-(\ref{eq:2.2}) restore the charging and discharging rates for the last period prior to the planning time $j$.
\begin{eqnarray}
 X_{c,i}^{\phi_j-1} &=LC_i^{\{j\}}, & \quad \forall i \in E^j \label{eq:2.1}\\
 X_{d,i}^{\phi_j-1} &=LD_i^{\{j\}}, & \quad \forall i \in E^j_{v2g} \label{eq:2.2}
\end{eqnarray}
\noindent The rest of the constraints in dynamic model are as follows:
\begin{eqnarray}
 & SOC_i^{\phi_j} = SOC_{init,i} , \quad \forall i \in E^j \label{eq:2.3}\\
 & SOC_i^t = SOC_i^{t-1} + \eta_i^{c} P_i^{c} X_{c,i}^{t-1} , \quad \forall i \in E^j_{g2v}, \quad \forall t \in \{\phi_j+1, ..., t_i^{dep}\} \label{eq:2.4}\\
 & SOC_i^t = SOC_i^{t-1} + \eta_i^{c} P_i^{c} X_{c,i}^{t-1} - P_i^{d} X_{d,i}^{t-1}/\eta_i^{d}, \quad \forall i \in E^j_{v2g}, \;\; \forall t \in \{\phi_j+1, ..., t_i^{dep}\} \label{eq:2.5}\\
 & SOC_{init,i} + P_i^{c}\eta_i^{c} (t_i^{dep}-\phi_j) \geq SOC_{desired,i} - MZ_i , \quad \forall i \in E^j \label{eq:2.6}\\
 & X_{c,i}^t \geq Z_i, \quad \forall i \in E^j ,\;\; \forall t \in T_i^p \label{eq:2.7}\\
 & SOC_i^{t_i^{dep}} \geq SOC_{desired,i} - M Z_i , \quad \forall i \in E^j \label{eq:2.8}\\
 & X_{c,i}^{t} \leq 1-Y_{c,i}^t , \quad \forall i \in E^j_{v2g}, \;\; \forall t \in T_i^p   \label{eq:2.9}\\
 & X_{d,i}^{t} \leq 1-Y_{d,i}^t , \quad \forall i \in E^j_{v2g}, \;\; \forall t \in T_i^p   \label{eq:2.10}\\
 & Y_{c,i}^t +Y_{d,i}^t = 1, \quad \forall i \in E^j_{v2g}, \;\; \forall t \in T_i^p \label{eq:2.11}\\
 & SOC_i^t \geq SOC_{min,i} , \quad \forall i \in E^j \setminus B, \;\; \forall t \in T_i^p \label{eq:2.12}\\
 & X_{c,i}^t = 1, \quad \forall i \in B, \;\; \forall t \in \{t_i^{arr},..., t_i^{arr} + T_{min,i}\} \label{eq:2.13}\\
 & SOC_i^t \geq SOC_{min,i}, \quad \forall i \in B, \;\; \forall t \in \{t_i^{arr} + T_{min,i} + 1,..., t_i^{dep}\} \label{eq:2.14}\\
 & SOC_i^t \leq SOC_{cap,i}, \quad \forall i \in E^j, \;\; \forall t \in \{\phi_j,..., t_i^{dep}\} \label{eq:2.15}\\
 & 0 \leq X_{c,i}^t \leq 1, \quad \forall i \in E^j, \;\; \forall t \in T_i^p \label{eq:2.16}\\
 & 0 \leq X_{d,i}^t \leq 1, \quad \forall i \in E^j_{v2g}, \;\; \forall t \in T_i^p \label{eq:2.17}\\
 & Y_{c,i}^1,\; Z_i \;\; Binary, \quad \forall i \in E^j, \;\; \forall t \in T_i^p\label{eq:2.18}\\ 
 & Y_{d,i}^t \;\; Binary, \quad \forall i \in E^j_{v2g}, \;\; \forall t \in T_i^p \label{eq:2.19}\\ 
 & \Omega^t,\; G^t \geq0, \quad \forall t \in \{\phi_j, ..., \tau_{max}^j-1\} \label{eq:2.20}
\end{eqnarray}

\section{Simulation} \label{sec:simulation}
We perform comprehensive simulations to examine the performance of the proposed controlled charge/discharge scheduling algorithm.

\subsection{Simulation Settings} \label{sec:simulationSetting}
In our simulation, a day is evenly divided into 96 time intervals of $\Delta t = 15$ minutes. For instance, $t=0$ denotes midnight and the start of scheduling period and $t=96$ is midnight of the next day. Thus, the decision periods are $T=\{0, 1, \dots, 95\}$. If a vehicle, for instance, arrives at time 6:00 am and departs at 12:15 pm, then $t^{arr} = 24$ and $t^{dep}=49$. The algorithm's goal is to decide the battery charge/discharge rate of EV $i$ during interval $t$ ($X_{c,i}^t,\;\; X_{d,i}^t$).

\noindent{\textit{Generation \& Consumption data:} The variations for electricity price and wind power generation are captured in ten different scenarios at various months of the year. Hourly wind power production is simulated using Grid Lab System Advisor Model (SAM) in northern California for a single wind turbine based on specifications of Endurance X33 turbine with 230 kW power capacity in ten consecutive days for all ten scenarios \cite{SAM}. Hourly electricity price is collected using historical Locational Marginal price (LMP) data for day-ahead market at node "PLAINFLD\_6\_N001" from California Independent System Operator (CAISO) for the same days \cite{CAISO}. The ten scenarios include five months in Spring and Summer seasons, and five months in Fall and Winter seasons. Thus, we consider a total of 100 days, each day has random arrival and departure behavior as well as different wind and electricity price profiles. Here, we propose that the discharged energy of an EV is sold with a price of 90\% of the real-time electricity price. Since the revenue made by the EV owners in discharge mode is equal to the charge cost of discharge energy for the EVs in charge mode, we did not include that in our mathematical model.

\noindent\textit{EV-related data:} The owners behaviors are also simulated for ten scenarios. Considering a total of 100 EV trips per day, the battery and charging characteristics of EVs are based on the specifications of 10 different EVs available in the market in 2018 \cite{EVs:2018}. The EV battery info is summarized in table \ref{tab:EVs}. EV arrival times are captured using data from the National Household Travel Survey (NHTS) \cite{NHTS:2017}. It is assumed that fifty of the EV charging occur at workplace and the other fifty are vehicles parking at home. Thus, we used travel time data related to trips with home or work purposes from NHTS. Figure \ref{fig:subfigure1} shows the probability of trip to home for each hour of the day. In Figure \ref{fig:subfigure2}, the x-axis represents hour of the day, and the probability of trip to workplace is on the y-axis. In each scenario of the simulation, the optimal schedule algorithm is run for 10 days. The EV plug-in period is modeled as discrete uniform distribution between 4 and 12 hours with a resolution of 15 minutes. The desired level of battery at departure is modeled using uniform distribution in the range of 0.75 to 0.95 of battery capacity. We assume charging and discharging efficiency of 90\%. The initial state of charge is also assumed to follow a uniform distribution between 0 and 0.65 of the battery capacity, while the minimum required state of charge ($SOC_{min}$) for all vehicles is assumed to be 5 kWh that in case of emergency seems sufficient to reach a distance of approximately 25 miles. For battery degradation cost, as in \cite{he2012optimal:10ev}, values of 0.05 and 0.1 $\cent$ are used for $\alpha$ and $\beta$, respectively.
    \begin{table}[h!]
    \caption{EV characteristics}
    \centering
    \begin{tabular}{  c  c  c  c  c }
    \hline
    \thead{Vehicle} & \thead{Acceptance\\ Rate\\ $AR$} & \thead{Battery\\ Size \\ $SOC_{cap}$} & \thead{Charger\\ Capacity\\ $CP$} & \thead{Battery\\ Cost\\ $C_{bat}$(\$)}\\
    \hline\hline
    BMW i3 2017 & 7.4 & 32 & 7.7 & 4,640\\
    Ford Focus EV & 6.6 & 23 & 7.7 & 3,500\\
    Ford Focus EV 2017 & 6.6 & 33.5 & 7.7 & 4,850\\
    Nissan Leaf S 2016 & 6.6 & 24 & 7.7 & 3,500 \\
    Nissan Leaf 2017 & 6.6 & 30 & 7.7 & 4,350 \\
    VW e-Golf 2017 & 7.2 & 35.8 & 7.7 & 5200\\
    Chevy Bolt & 7.2 & 60 & 7.7 & 8,700\\
    Tesla Model S 70 Single & 9.6 & 70 & 11.5 & 10,150\\
    Tesla Model X 75 Dual & 17.2 & 75 & 15.4 & 10,900\\
    Tesla Model S 90 Dual & 19.2 & 90 & 15.4 & 13,000\\
    \hline
    \end{tabular}
    \label{tab:EVs}
    \end{table}
\begin{figure}[t!]
\centering
\caption{Distribution of arrival times at a) home b) workplace}
\begin{subfigure}{.475\textwidth}
    \centering
    \caption{Distribution of arrivals at home}
    \includegraphics[width=0.9\linewidth]{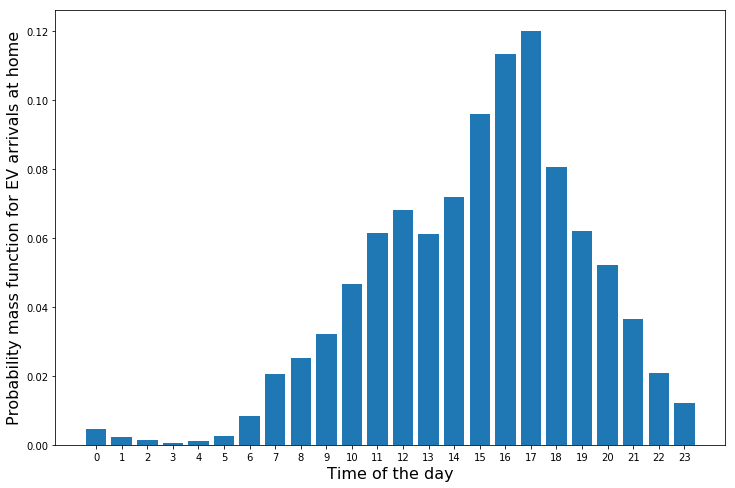}
    \label{fig:subfigure1}
\end{subfigure}%
\begin{subfigure}{.475\textwidth}
    \centering
    \caption{Distribution of arrivals at workplace}
    \includegraphics[width=0.9\linewidth]{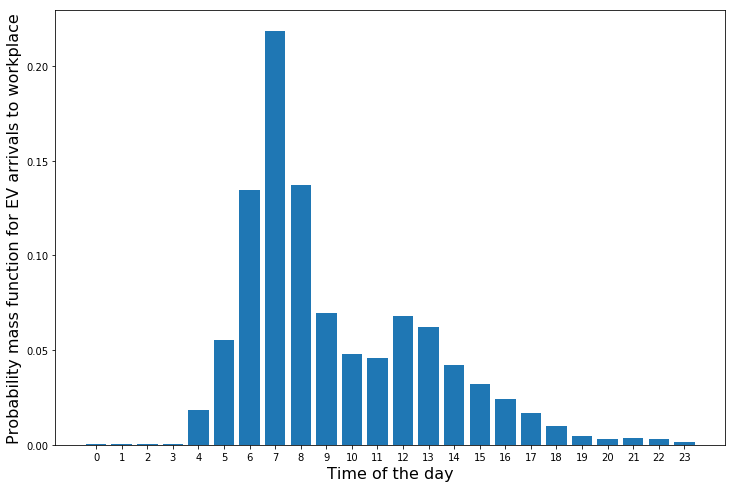}
    \label{fig:subfigure2}
\end{subfigure}
\label{fig:probMassFunction}
\end{figure}
\subsection{Results and Performance Analysis}
To evaluate the performance of the proposed approach, the results of the dynamic charging algorithm are compared with the static and BAU charging scenario. All models are run for the ten scenarios mentioned in Subsection \ref{sec:simulationSetting}. The simulation is coded in Python 3.7 using Gurobi optimizer and the results for optimal solutions are shown in Figure \ref{fig:ComparisonAll}. Comparing the results, the proposed charging algorithm causes significant improvement in all the objective measures including energy supplied from the grid, wind utilization, and total cost for EV owners. The total cost consists of charging cost, degradation cost, and revenue made from selling discharging energy.

\begin{figure} [h!]
    \centering
    \caption{\small Performance evaluation for all three charging cases, a) total grid supply, b) wind utilization, and c) total cost} 
    \begin{subfigure}[b]{0.475\textwidth}
        \centering
        \caption[]%
        {{\footnotesize Total grid supply (kWh)}}
        \includegraphics[width=\textwidth]{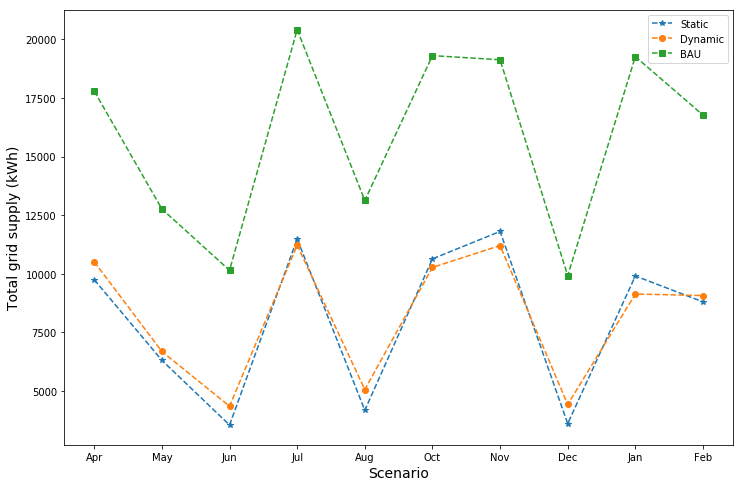}
        \label{fig:GS_all}
    \end{subfigure}
    \begin{subfigure}[b]{0.465\textwidth}  
        \centering 
        \caption{{\footnotesize Wind utilization (\%)}}    
        \includegraphics[width=\textwidth]{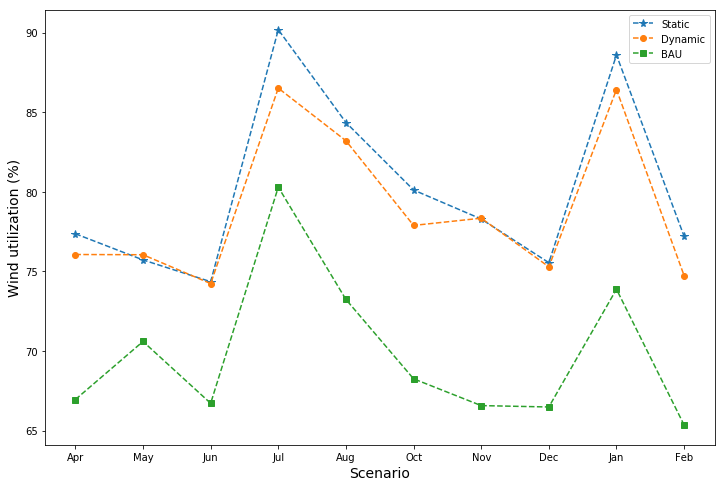}
        \label{fig:WU_all}
    \end{subfigure}
    \begin{subfigure}[b]{0.475\textwidth}   
        \centering 
        \caption[]%
        {{\footnotesize Total cost ($\cent$)}} 
        \includegraphics[width=\textwidth]{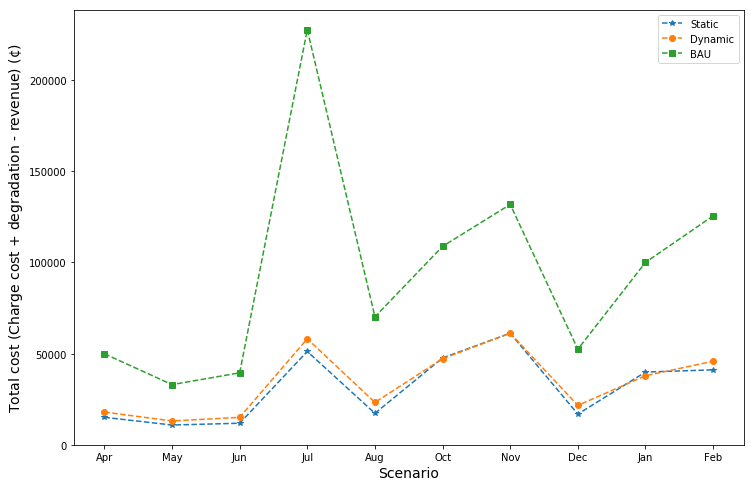}
        \label{fig:TC_all}
    \end{subfigure}
    \label{fig:ComparisonAll}
\end{figure}

The results demonstrates that the dynamic model achieves a similar performance compared to the static scheduling model. It is worthwhile to mention that the performance of the static case for some objective metrics is worse than the dynamic case for a few scenarios due to the fact that we consider a multi-objective optimization. However, the total objective value for the static case is always better than the dynamic case.

To validate our degradation model, Figure \ref{fig:ComparisonDegC} shows that considering the linear degradation model described in \ref{battDeg} leads to very similar results compared to the quadratic model.

\begin{figure}[t!]
\centering
\caption{Comparison of the models with quadratic vs linear degradation function.}
\begin{subfigure}{.475\textwidth}
    \centering
    \caption{Total cost}
    \includegraphics[width=\linewidth]{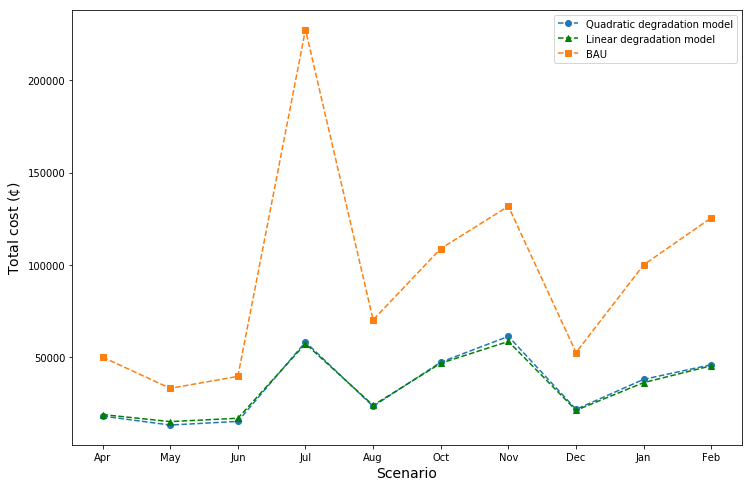}
    \label{fig:subfigure11}
\end{subfigure}%
\begin{subfigure}{.475\textwidth}
    \centering
    \caption{Wind utilization}
    \includegraphics[width=\linewidth]{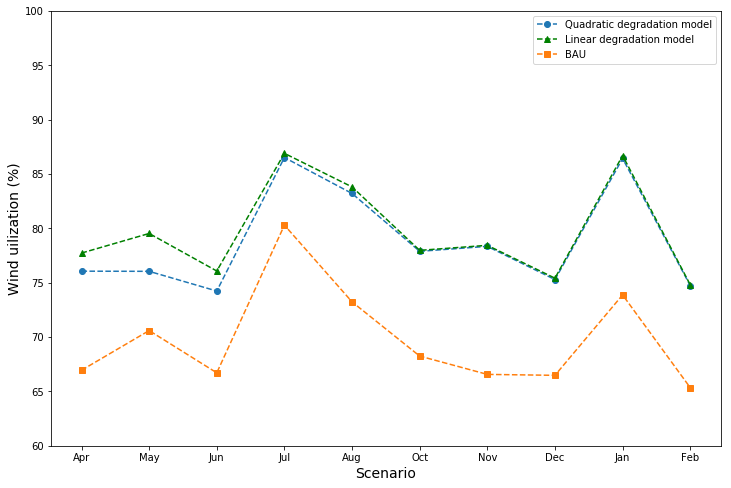}
    \label{fig:subfigure22}
\end{subfigure}
\label{fig:ComparisonDegC}
\end{figure}

\subsubsection{Discussion of Objective Function Hyper-Parameters}
Since we consider a weighted multi-objective optimization problem, we need to evaluate our objective for different values of hyper-parameters. The values of $\delta$ and $\lambda$ are parameters that determine the importance of the corresponding term in the objective function. We will analyze the value of these two hyper-parameters for our objective measures.

\noindent \textit{Performance evaluation under different values of $\delta$:} Giving a constant weight for wind curtailment does not capture the dynamic of our model well. Thus, we consider the penalty relative to the electricity price such that wind curtailment penalty and grid supply have same ratio through all periods. Here, $\delta$ is considered as multiplier of $pr^t$ in the third term of objective function. A small value of $\delta$ gives little weight to wind curtailment minimization, while making sure that grid supply is as low as possible. If a large value is given to $\delta$, wind curtailment is penalized more causing to increase wind utilization in the cost of higher discharged energy and higher degradation cost, which as a result, leads to higher total charging cost. Figure \ref{fig:sigmaEvaluation} shows the wind utilization percentage (a) and total cost (b) for multiple values of $\delta$ ($\lambda$ value is fixed and equal to 1). If wind utilization is of top priority, then a high value should be used for $\delta$, and vice versa. Experimentally, a value of 0.25 seems to have the best performance among all options.

\begin{figure}
\centering
\caption{Performance evaluation under different values of $\delta$}
\begin{subfigure}{.5\textwidth}
    \centering
    \caption{Wind utilization.}
    \includegraphics[width=\linewidth]{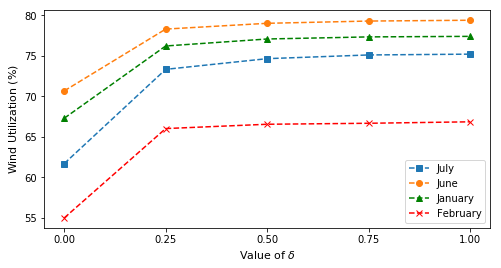}
    \label{fig:delta_wind}
\end{subfigure}%
\begin{subfigure}{.5\textwidth}
    \centering
    \caption{Total cost relative to the $\delta =0$ case.}
    \includegraphics[width=\linewidth]{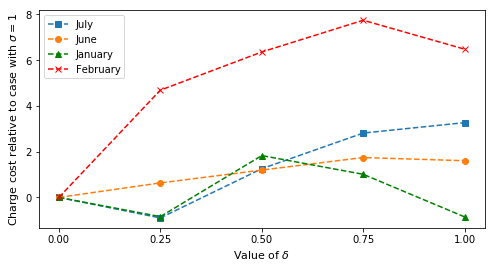}
    \label{fig:sigmaCost}
\end{subfigure}
\label{fig:sigmaEvaluation}
\end{figure}

\noindent \textit{Performance evaluation under different values of $\lambda$:} The value of $\lambda$ indicates the tolerance of EV owners for their battery degradation. $\lambda = 0$ means high tolerance, while value of 1 means battery degradation cost is as important as charging cost for the owner. Figure \ref{fig:lambdaEvaluation} plots the battery degradation cost, total cost, and wind utilization percentage for multiple values of $\lambda$. As seen in these plots, a low value leads to higher wind utilization in cost of more degradation and total cost.

\begin{figure}[t!]
\centering
\caption{Performance evaluation under different values of $\lambda$}
\begin{subfigure}{.325\textwidth}
    \centering
    \caption{Wind utilization}
    \includegraphics[width=\linewidth]{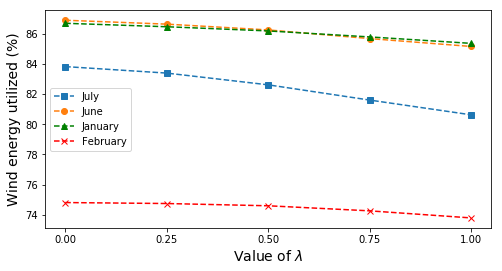}
    \label{fig:sigma_wind}
\end{subfigure}%
\begin{subfigure}{.325\textwidth}
    \centering
    \caption{Degradation cost}
    \includegraphics[width=\linewidth]{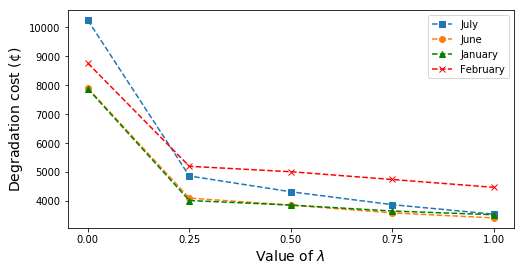}
    \label{fig:deltaCost}
\end{subfigure}
\begin{subfigure}{.325\textwidth}
    \centering
    \caption{Total cost}
    \includegraphics[width=\linewidth]{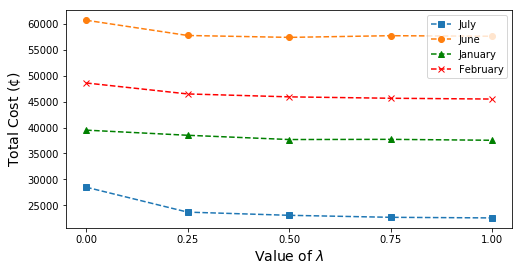}
    \label{fig:lambdaTC}
\end{subfigure}
\label{fig:lambdaEvaluation}
\end{figure}

\subsubsection{V2G Benefits}

To assess the benefits of V2G in this MG, we compare the scenarios in which the percentage of vehicles participating in V2G versus G2V varies. Different values of $R_{v2g}$ is considered to evaluate to the benefits of V2G to both the owner and the grid. The performance evaluation under different value of $R_{v2g}$ is shown in Figure \ref{fig:RV2G}. In order to maintain clarity, the plots are shown for four months instead of all ten scenarios. In case where all vehicles participate in V2G, the grid benefits from increasing the wind utilization and reducing the total energy curtailment. From the owners perspective, the degradation cost increases as more vehicles discharge their energy; however, the total cost decreases mainly due to the fact that the energy purchased from discharged energy reduces the need for expensive supply from the external grid. In V2G case, the EVs store energy from wind source and in low wind periods, they inject energy back to the system to charge other EVs. The breakdown of total cost for the EV owners from the July scenario is provided in Figure \ref{fig:COSTREvenue_pvalue}. The charge cost consists of the cost of energy purchased from the grid and the discharge energy. The revenue and the discharged energy cost cancel each other out; however, since the total supply from the grid decreases, the total cost decreases as well. 

\begin{figure} [t!]
    \centering
    \caption{\small Performance evaluation under different values of $R_{V2G}$ a) wind utilization, b) total cost, c) degradation cost, d) energy discharged } 
    \begin{subfigure}[b]{0.475\textwidth}
        \centering
        \vspace{-0.5in}
        \caption[Network2]%
        {{\footnotesize Wind utilization (\%)}} 
        \includegraphics[width=\textwidth]{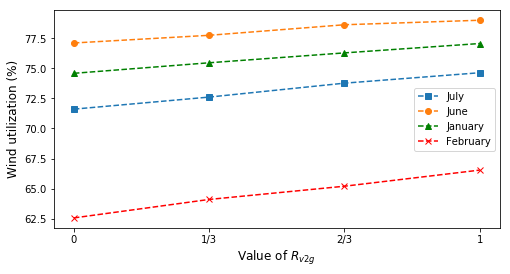}
        \label{fig:WindUtilPvalue}
    \end{subfigure}
    \begin{subfigure}[b]{0.475\textwidth}  
        \centering 
        \caption[]%
        {{\footnotesize Total cost ($\cent$)}} 
        \includegraphics[width=\textwidth]{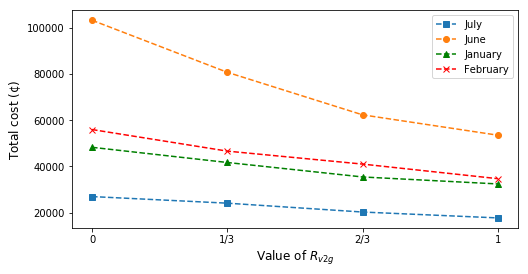}
        \label{fig:TotalCostPvalue}
    \end{subfigure}
    \begin{subfigure}[b]{0.475\textwidth}   
        \centering 
        \caption[]%
        {{\footnotesize Total degradation cost ($\cent$)}} 
        \includegraphics[width=\textwidth]{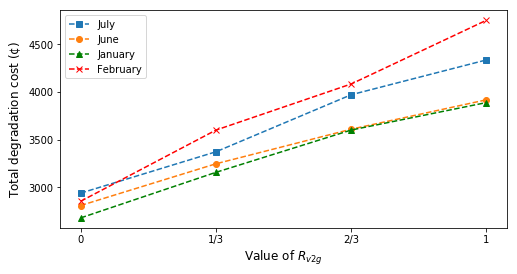}
        \label{fig:DegCost}
    \end{subfigure}
    \begin{subfigure}[b]{0.475\textwidth}   
        \centering 
        \caption[]%
        {{\footnotesize Total energy discharged (kWh)}} 
        \includegraphics[width=\textwidth]{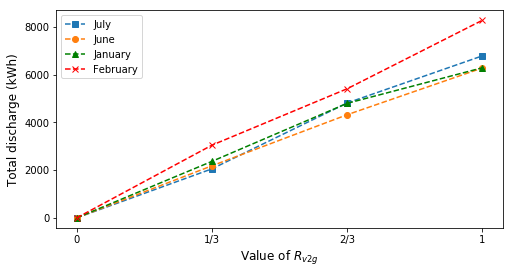}
        \label{fig:energyDCH}
    \end{subfigure}
    \label{fig:RV2G}
\end{figure}

\begin{figure}[t]
    \centering
    \caption{Charge cost and revenue from discharge for different value of $R_{V2G}$ (July scenario)}
    \centerline{\includegraphics[scale=0.75]{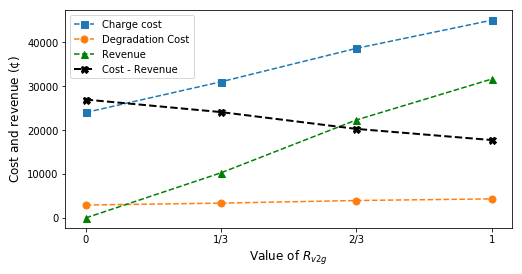}}
    \label{fig:COSTREvenue_pvalue}
\end{figure}

\subsubsection{Wind Forecast Uncertainty} \label{sec:intermittentWind}
In reality, the wind forecast for the following couple of hours is not perfect. Thus, we consider a model that with updated wind forecast, the forecast is only accurate for the current $\Delta j$ planning interval. We believe that in the dynamic model, since the scheduling of EVs are updated every one hour, the effect of intermittent wind generation is significantly reduced and the model can accommodate the uncertainty of wind production to a higher degree. To prove that, we consider a model with perfect wind forecast for the next one hour, beyond that, we use a discrete time Markov Decision Process (MDP) model to estimate the wind generation \cite{kwon2015meeting:34'}. In a simple case, we discretize the state space for wind scenario and consider a MDP with 20 discrete states and estimate the transition probability based on historical (training) data. The forecast for the next $k$ hours is estimated by  
\begin{equation}
    W_{forecast}^{t+k} = \displaystyle \sum_{w \in S_w} P_{W^t, w}^k \times w, \quad k = \{1,2,3, ...\}
    \label{windMDP}
\end{equation}
where, the state space is denoted by $S^w$ and consists of 20 wind scenario representatives. $P_{W^t, w}^k$ is the transition probability from current state of wind generation to state $w$ in $k$ transition steps. The first 15 days of each month are used as training to calculate the transition probability matrix, and the problem is optimized for the next 10 days with MDP wind forecast. The results are shown in Figure \ref{fig:ComparisonStoch}, where it is observed that the dynamic model with imperfect wind forecast performs similar to the perfect forecast scenario. Also, note that we only considered a very simple MDP process that might not be an accurate forecast. A better forecast for wind generation is likely to improve the quality of solutions. Same procedure can also be applied to electricity price to predict real-time price, which is omitted from this paper due to page limit restrictions. 

\begin{figure*}[h]
    \centering
    \caption[ ]
    {\small Comparison of dynamic model with perfect forecast vs MDP wind forecast.}
    \begin{subfigure}[b]{0.325\textwidth}
        \centering
        \caption[]%
        {{\footnotesize Total grid supply (kWh)}}    
        \includegraphics[width=\textwidth]{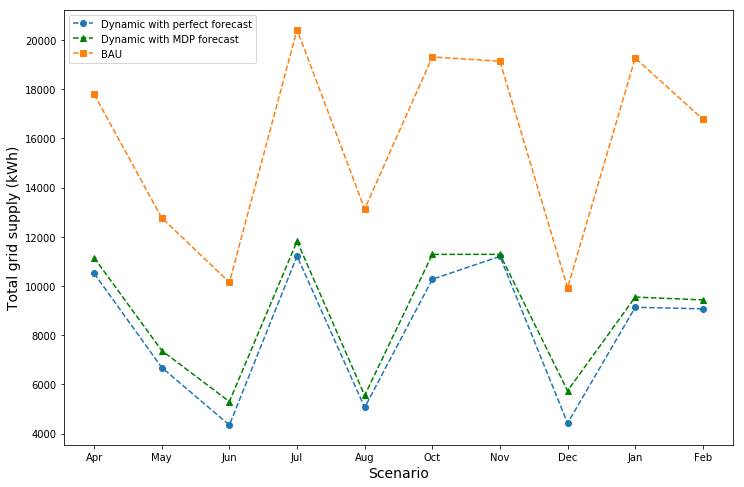}
        \label{fig:GS_stoch}
    \end{subfigure}
    \begin{subfigure}[b]{0.315\textwidth}  
        \centering 
        \caption[]%
        {{\footnotesize Wind utilization (\%)}}   
        \includegraphics[width=\textwidth]{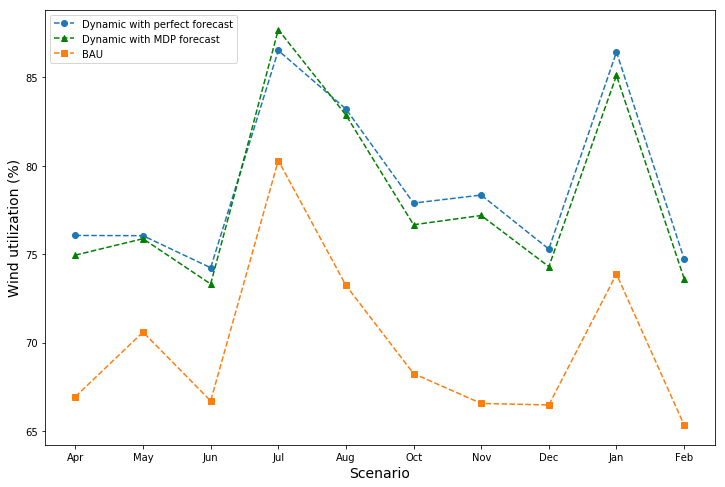}
        \label{fig:WU_stoch}
    \end{subfigure}
    \begin{subfigure}[b]{0.325\textwidth}   
        \centering 
        \caption[]%
        {{\footnotesize Total cost ($\cent$)}}    
        \includegraphics[width=\textwidth]{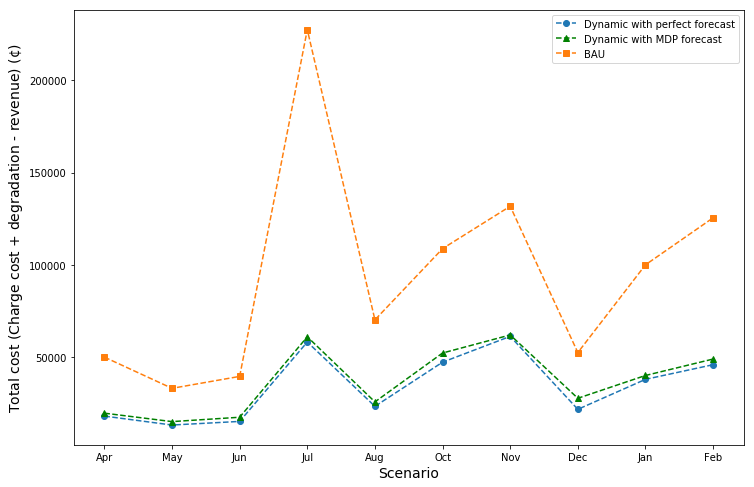}
        \label{fig:TC_stoch}
    \end{subfigure}
    \label{fig:ComparisonStoch}
\end{figure*}

\section{Conclusion} \label{conclusion}
In this paper, an optimization algorithm for EV charge/discharge scheduling is proposed to support wind energy integration using a rolling horizon optimization method. The problem is formulated as a MIQP and the results show significant improvement in our proposed approach compared to the BAU case in terms of total charging cost, wind utilization, and demand from conventional generators. The main contribution here is to design the scheduling algorithm that aggregators can exploit the presence of advanced communication technology in smart grid, flexibility of EV drivers, and the V2G technology to support high integration of intermittent wind energy into the power system. Unlike some of the previous research, the proposed algorithm does not attempt to start or stop charging EVs frequently, which leads to battery degradation. The approach here is based on optimal scheduling of a large number of EVs depending on the availability of EVs and preferences, needs and flexibility of their owners. To mitigate the barriers in people participation in V2G, the battery degradation, minimum required level of charge, and/or financial incentives for the EV drivers are considered. Furthermore, in this work, multi-objective optimization is considered to maximize wind utilization, minimize the demand from conventional generators, and minimizing charging cost while satisfying the driver needs and preferences. A simulation of the proposed algorithms for different scenarios of EV characteristics, arrivals, departures, and charging requirements are performed to check the quality of solutions and schedules. The results show that the proposed model leads to significant improvement in all metrics and benefits both the owner and the grid. Moreover, the results indicate that frequent updates of available wind power (and electricity price) in the deterministic problem significantly reduce the effects of forecast uncertainty. Future research is needed to consider a large-scale grid with multiple generators and consumers with the presence of battery storage and demand side management for residential load to evaluate the users flexibility in power system. A financial incentive framework should also be developed to encourage more people to participate in V2G.


\end{document}